\documentclass[12pt]{amsart}
\usepackage{amssymb,bm,fullpage,graphicx,hhline,longtable,mathtools,mathrsfs,tikz,tikz-cd,upref,xurl}
\usepackage{anyfontsize}
\usepackage{newtxtext,newtxmath}
\usepackage{hyperref}

\newtheorem{theorem}{Theorem}[section]

\theoremstyle{definition}

\theoremstyle{remark} 
\newtheorem*{remark}{Remark}

\numberwithin{equation}{section}

\let\originalleft\left
\let\originalright\right
\renewcommand{\left}{\mathopen{}\mathclose\bgroup\originalleft}
\renewcommand{\right}{\aftergroup\egroup\originalright}

\newcommand{\p}[2]{{#1}^{({#2})}}

\makeatletter
\@namedef{subjclassname@2020}{%
  \textup{2020} Mathematics Subject Classification}
\makeatother

\begin{document}

\title[Parity results  concerning the generalized divisor function]{Parity results concerning the generalized divisor function involving small prime factors of integers}

\author{Krishnaswami Alladi}

\author{Ankush Goswami}

\address{Department of Mathematics, University of Florida, Gainesville, FL, 32611, USA}
\email{alladik@ufl.edu}

\address{School of Mathematical and Statistical Sciences, the University of Texas Rio Grande Valley, 1201 W. University Dr., Edinburg, TX, 78539}
\email{ankush.goswami@utrgv.edu, ankushgoswami3@gmail.com}



\dedicatory{To the memory of Professor M. V. Subbarao on the occasion of his Centenary}
\maketitle

\begin{abstract}
  Let $\nu_y(n)$ denote the number of distinct prime factors of $n$ that are
  $<y$. For $k$ a positive integer, and for $k+2\leq y\leq x$, let $S_{-k}(x,y)$
denote the sum
\begin{eqnarray*}
S_{-k}(x,y):=\sum_{n\leq x}(-k)^{\nu_y(n)}.
\end{eqnarray*}
In this paper, we describe our recent results on the asymptotic
behavior of $S_{-k}(x,y)$ for $k+2\leq y\leq x$, and $x$ sufficiently large. There is a crucial difference in the asymptotic behavior of $S_{-k}(x,y)$
when $k+1$ is a prime and $k+1$ is composite, and this makes the problem
particularly interesting. The results are derived utilizing a combination of
the Buchstab-de Bruijn recurrence, the Perron contour integral method, and
certain difference-differential equations. We present a summary of our results against the background of earlier work of the first author on
sums of the M\"{o}bius function over integers with restricted prime factors
and on a multiplicative generalization of the sieve.  
\end{abstract}

\section{Introduction and Background}
\noindent Let $\gamma_y(n)$ (resp. $\nu_y(n)$) denote the product of the (resp. the number of) distinct prime factors $<y$ of $n$. The count of the number of distinct positive divisors of $\gamma_y(n)$, that is, the number of ordered positive integer solutions $(a,b)$ to the equation $ab=\gamma_y(n)$ is $2^{\nu_y(n)}$. Similarly, for $k>2$, $k^{\nu_y(n)}$ counts the number of ordered $k$-tuples, $(d_1, d_2,\cdots, d_{k-1}, d_k)$ of positive divisors of $\gamma_y(n)$ such that $d_1\cdot d_2\cdots d_{k-1}\cdot d_k=\gamma_y(n)$.
Hence $(-k)^{\nu_y(n)}$ counts these tuples with a sign (positive or negative) depending on the parity (even or odd) of $\nu_y(n)$. In this paper, we discuss the asymptotic behavior of the sum
\begin{eqnarray}\label{smain}
S_{-k}(x,y):=\sum_{n\leq x}(-k)^{\nu_y(n)},
\end{eqnarray}
uniformly for $k+2\leq y<x$.

There are striking differences in the asymptotic behavior of $S_{-k}(x,y)$ depending on $k=p-1$ (for some prime $p$) and $k\neq p-1$ (for any prime $p$) which is what prompted this study. The lower bound $k+2$ for $y$ is to ensure that
all primes $p\leq k+1$ are accounted for in the study of $\nu_y(n)$. The analysis of $S_{-k}(x,y)$ involves a careful study of the asymptotic behavior separately for ``small $y$'' and ``large $y$''. The terms ``small'' and ``large'' used in a technical sense will be made precise later. For small $y$
we use the Perron integral to establish the asymptotic behavior of $S_{-k}(x,y)$, and for large $y$, we employ the Buchstab-de Bruijn iteration method to asymptotically estimate $S_{-k}(x,y)$ by induction on $\lfloor\alpha\rfloor$ where $\alpha=\log x/\log y$. It turns out that the two regions for small and large $y$ 
overlap and thus we can determine the asymptotic behavior of $S_{-k}(x,y)$ uniformly for $k+2\leq y\leq x$ by comparing the estimates in these regions.

When $k=p-1$, the main term in the asymptotic estimate of $S_{-k}(x,y)$ for small $y$ vanishes and we only get an upper bound using Perron's method in this region, but we do get a uniform asymptotic estimate of $S_{-k}(x,y)$ for large $y$
using the de Bruijn-Buchstab iteration. When $k\neq p-1$, the main term does not vanish for small $y$ and so we obtain asymptotic estimates for $k+2\leq y\leq x$. To understand these
intricacies in the asymptotic behavior of $S_{-k}(x,y)$, we need to describe certain earlier results as a background.

\subsection{Sums of the M\"{o}bius function} The fundamental tool used to study the parity of the number of prime factors $\nu(n)$ on the square-free numbers is the M\"{o}bius function $\mu(n)$ which is defined as $\mu(n)=(-1)^{\nu(n)}$, when $n$ is square-free, and 0 otherwise. Bounds for $M(x)$, the sum of $\mu(n)$ for $1\le n\le x$,  are closely tied to the error term in the strong form of the Prime Number Theorem (PNT), and to the Riemann Hypothesis as well (see Tenenbaum \cite{GT}). Motivated by this strong link, Alladi \cite{A1}, \cite{A2} studied the asymptotic behavior of the following two sums:
$$
M(x,y)=\sum_{n\le x, p(n)>y}\mu(n),
$$
where $p(n)$ is the smallest prime factor of $n$ if $n>1$, and
$p(1)=\infty$, and
$$
M^*(x,y)=\sum_{n\le x, P(n)<y}\mu(n),
$$
where $P(n)$ is the largest prime factor on $n$ if $n>1$ and $P(1)=1$. Note that $M(x,1)=M^*(x,x)=M(x)$. The functions $M(x,y)$ and $M^*(x,y)$ are weighted versions of the well known functions $\Phi(x,y)$ - the number of uncancelled elements in the sieve of Eratosthenes, and $\Psi(x,y)$ - the number of integers up to $x$ free of prime factors $>y$. 

In \cite{A1} it is shown that if $\alpha:=\log x/\log y>1$, is fixed, then
$$
M(x,y)=\frac{x\cdot m(\alpha)}{\log y} + O\left(\frac{x}{\log^2y}\right),
$$
where $m(\alpha)$ satisfies a difference-differential equation and is the derivative of the famous Dickman function $\rho(\alpha)$ that occurs in the asymptotic estimate of $\Psi(x,y)$. Similarly in \cite{A2} it is shown that
\begin{equation}\label{M*}
M^*(x,y)=\frac{x\cdot f(\beta)}{\log y}+\frac{x\cdot m^*(\alpha)}{\log^2y}+
O\left(\frac{x}{\log^3 y}\right),    
\end{equation}
where $\beta=x/y$, and $m^*(\alpha)$ is the derivative of the Buchstab function $w(\alpha)$ that occurs in the asymptotic estimate of $\Phi(x,y)$. It turns out that $f(\beta)\to 0$ rapidly as $\beta\to\infty$, and so the striking feature about $M^*(x,y)$ is that when $\beta>1$ is ``finite'', that is when $y$ is like a constant times $x$, $M^*(x,y)$ behaves like $x/\log y$, but when $\beta$ is large, which happens if $\alpha>1$ is fixed, then $M^*(x,y)$ behaves like $x/\log^2 y$. It is this type of change in asymptotic behavior, that happens for $S_{-k}(x,y)$, but in a much more intricate manner, as will be seen below. We also point out that in \cite{A1} and \cite{A2}, uniform asymptotic estimates for long ranges of $\alpha$ are established for $M(x,y)$ and $M^*(x,y)$. See de Bruijn \cite{deB1} and \cite{deB2} for initial results related to the uniform asymptotic estimates of $\Phi(x,y)$ and $\Psi(x,y)$ respectively. A comprehensive survey of subsequent results on $\Psi(x,y)$ and related functions is given in Hildebrand and Tenenbaum \cite{H-T}. 

\subsection{A multiplicative generalization of the sieve} Before discussing $S_{-k}(x,y)$ for all positive integers $k$, we will state the results in Dhavakodi's 1992 PhD thesis \cite{Dha} for $S_{-1}(x,y)$, but to motivate this, we will first discuss a multiplicative generalization of the sieve due to the first author \cite{A3}.
  
Let $\mathscr{A}\subseteq\mathbb{N}$ be a subset of the positive integers, and let $\mathscr{A}(x):=\mathscr{A}\cap[1,x]$. Let $\mathscr{P}$ be a set of primes and for $y>0$ , define
$$
P_y:=\prod_{\substack{p\in\mathscr{P}\\p<y}}p.
$$
The classical sieve problem is to estimate the quantity 
\begin{eqnarray}
\mathscr{S}(\mathscr{A};\mathscr{P},x,y):=\sum_{a\in\mathscr{A},a\le x, (a,P_y)=1}1.
\end{eqnarray}
To suitably rewrite (1.3) for estimation, one considers
$$
\mathscr{A}_d=\mathscr{A}_d(x):=\{a\in\mathscr{A}: a\le x, a\equiv 0 \pmod{d}\}
$$
so that
$$
\mathscr{S}(\mathscr{A};\mathscr{P},x,y)
=\sum_{a\in\mathscr{A},a\le x}\sum_{d|(a,P_y)}\mu(d)=\sum_{d|P_y}\mu(d)|\mathscr{A}_d|.
$$
Estimates are usually obtained (see \cite{H-R}) under suitable assumptions on
$\mathscr{A}_d$, such as:
\begin{equation*}
|\mathscr{A}_d|=\dfrac{X\omega(d)}{d}+R_d,
\end{equation*}
where $X$ is an approximation to $|\mathscr{A}(x)|$, $\omega$ is multiplicative, and the remainder $R_d$ can be controlled
in an average sense. Since the divisors of $P_y$ are too numerous when $y$ is large, the contribution due to the sum over $R_d$  cannot be controlled. Thus
Brun had the brilliant idea to consider upper and lower ``sifting functions''
$\chi^+$ and $\chi^-$ that will satisfy
\begin{equation}
\sum_{d|n}\mu(d)\chi^-(d)\leq\sum_{d|n}\mu(d)\leq\sum_{d|n}\mu(d)\chi^+(d),
\end{equation}
for all $n$. The $\chi^+$ and $\chi^-$ are chosen to vanish often enough to
keep the remainder term in check, and at the same time not affect the size
of the main term appreciably. Using Brun's sieve, it can be shown that the
following estimate holds:
\begin{equation}
  \mathscr{S}(\mathscr{A};\mathscr{P}, x ,y)=X\prod_{\substack{p\in\mathscr{P}\\p<y}}\left(1-\dfrac{\omega(p)}{p}\right)\left(1+O(e^{-\alpha})\right)
  +O_R\left(\frac{X}{\log^R X}\right), \quad \text{for all} \quad R>0,
\end{equation}
where 
\begin{equation*}
\alpha=\dfrac{\log x}{\log y}\iff y=x^{1/\alpha}.
\end{equation*}
Thus, if $\alpha\rightarrow\infty$ with $x$, then this yields an asymptotic estimate for $\mathscr{S}(\mathscr{A};\mathscr{P}, x, y)$.

It was such an asymptotic estimate for large $\alpha$ using Brun's sieve that was used in the proof of the celebrated Erd\H{o}s-Kac theorem. Even though sieve asymptotic estimates were employed in Probabilistic Number Theory, mainly the distribution of additive functions in the set of all positive integers $\mathbb{Z}^+$ was considered.

In order to study the distribution of additive functions in subsets, 
the first author considered the following sum \cite{A3}:
\begin{equation}\label{Sg}
S_g(\mathscr{A}(x),y):=\sum_{n\leq x,\;n\in\mathscr{A}}g_y(n)    
\end{equation}
where $g(n):=\prod_{p|n}g(p)$ is a strongly multiplicative function and
\begin{equation*}
  g_y(n):=\prod_{p|n,\;p<y}g(p)
\end{equation*}
is the \textit{trunction of $g$} at $y$. The sum in \eqref{Sg} can be viewed
as a multiplicative generalization of the sieve because the classical sieve problem is the special case when $g(p)=0$ for $p\in \mathscr P$ and $p<y$, with
$g(p)=1$, otherwise. The dual $g^*$ of $g$ is the multiplicative function
defined by
\begin{equation*}
g^*(n)=\prod_{p|n}g^*(p),\hspace{0.5cm}\text{where}\hspace{0.5cm}g^*(p)=1-g(p)    
\end{equation*}
for all $n$. The function $g^*$ occurs in the Buchstab identity for $S_g(\mathscr{A}(x),y)$, namely
$$
S_g(\mathscr{A}(x),y)=S_g(\mathscr{A}(x),y_1)+\sum_{y\leq p\ <y_1}g^*(p)S_g(\mathscr{A}_p(x),p).
$$
The functions $g$ and $g^*$ are related via the M\"{o}bius function relation  
\begin{equation*}
g^*(n)=\sum_{d|n}\mu(d)g(d)\iff g(n)=\sum_{d|n}\mu(d)g^*(d).    
\end{equation*}
Also
\begin{equation*}
g_y(n)=\sum_{d|(n,P_y)}\mu(d)g^*(d).  
\end{equation*}
Thus
$$\mathscr{S}(\mathscr{A};\mathscr{P}, x, y)=\sum_{n\in\mathscr{A}, n\leq x}
\sum_{d|(n,P_y)}\mu(d)g^*(d)=\sum_{d|P_y}\mu(d)g^*(d)|\mathscr{A}_d|.
$$
Of particular interest is the case $0\leq g\leq 1$ which means $0\leq g^*\leq 1$
as well and vice-versa. For such $g,\;g^*$, the first author noticed a {\it{monotonicity principle}}, namely that the Brun inequalities for $\chi^+$ and $\chi^-$
now imply that for all positive integers $n$,
\begin{equation}\label{genineq}
\sum_{d|n}\mu(d)g^*(d)\chi^-(d)\leq\sum_{d|n}\mu(d)g^*(d)\leq\sum_{d|n}\mu(d)g^*(d)\chi^+(d).
\end{equation}
What this means is that Brun's sieve can be used to estimate
$S_g(\mathscr{A}(x),y)$. More specifically, for all strongly multiplicative functions $g$ satisfying $0\leq g\leq 1$, it was shown in \cite{A3} that
$$
S_g(\mathscr{A}(x),y)= X\prod_{p<y}\left(1-\dfrac{\omega(p)g^*(p)}{p}\right)(1+O(e^{-\alpha}))+O_R\left(\frac{X}{\log^R X}\right), \quad \text{for all} \quad R>0
$$
holds uniformly. Another very unusual phenomenon was observed in \cite{A3}, namely for the {\it{Brun pure sifting functions}}, 
\begin{equation*}
\chi^+(d):=\begin{cases}
1,&\nu(d)\leq 2s\\
0,&\text{otherwise},
\end{cases}
\quad\quad\quad\quad
\chi^-(d):=\begin{cases}
1,&\nu(d)\leq 2s-1\\
0,&\text{otherwise},
\end{cases}
\end{equation*}
the {\it{monotonicity principle}} and the general inequalities \eqref{genineq} hold for all strongly multiplicative $g$ satisfying $0\leq g^*\leq 2$, that is for all $-1\leq g\leq 1$. In particular, the {\it{Brun pure sieve}} can be used to estimate $S_g(x,y)$ ($=S_g(\mathscr{A}(x),y)$ with $\mathscr{A}=\mathbb{N}$) when
$$
g(n)=(-1)^{\nu(n)}
$$
in the Brun pure sieve range
$$
\alpha > \log\log X.
$$
As far as we know this is the first instance when the sieve has been used to estimate sums of multiplicative functions which even take negative integer values. So this motivated the first author to consider estimates for $S_{-1}(x,y)$ which measures the parity of $\nu_y(n)$. In \cite{A3}, suitable bounds for $S_{-1}(x,y)$ were established. Subsequently, these bounds were
refined by Dhavakodi \cite{Dha} in his PhD thesis using analytic techniques. We now describe Dhavakodi's results which will set the stage for our work on $S_{-k}(x,y)$.

\subsection{Estimates for $S_{-1}(x,y)$} Dhavakodi \cite{Dha} used the de Bruijn-Buchstab method and so his starting estimate is in the range $\sqrt x < y\leq x$ which is  
$$
S_{-1}(x,y)=\frac{xh_{1}(\beta)}{\log y} +\frac{xh_2(\beta)}{\log^2y}+
  \frac{x\cdot m_{-1}(\alpha)}{\log^3y}+O\left(\frac{x}{\log^4y}\right),
$$
  where $\beta=x/y$ and $\alpha=\log x/\log y$. The functions $h_1$ and $h_2$
  rapidly tend to 0 as $\beta\to\infty$. Thus if $\alpha >1$ is fixed,
  then $S_{-1}(x,y)$ is of order of magnitude $x/\log^3y$. Compare this with the estimate for
  $M^*(x,y)$ in \eqref{M*}. 

  So with this starting estimate, he uses induction on $\lfloor\alpha\rfloor$ to establish the following estimate in the range $x^{1/3}<y<\sqrt{x}$:
  $$
  S_{-1}(x,y)=\frac{xh_{3}(\gamma)}{\log^2y}+
  \frac{x\cdot m_{-1}(\alpha)}{\log^3y}+O\left(\frac{x}{\log^4y}\right),
$$
  where $\gamma=x/y^2$. The main difference between $M^*(x,y)$ and $S_{-1}(x,y)$ is that for $S_{-1}(x,y)$, there are terms of size $x/\log y$, and $x/\log^2y$
  in the range $1<\alpha<2$ but these will shrink as $\alpha$ moves away from
  1, so that the dominant term is $x/\log^3y$. Next for $\alpha\geq 3$,
  there is no $x/\log y$ term, but there is the term of size $x/\log^2 y$.
  But this shrinks to 0 as $y$ moves away from $\sqrt{x}$ leaving $x/\log^3y$ as
  the dominating term. For $y<x^{1/3}$, namely for $\alpha >3$, the expansion
  starts with $x/\log^3y$ which is the dominating term. It is such a phenomenon that is seen more generally
  for $S_{-k}(x,y)$ as our results will show. 

  \section{The sum $S_z(x,y)$}

  For two complex valued functions $f$ and $g$, the Landau notation $f(x)=O(g(x))$ means that $|f(x)|<M|g(x)|$ for $x\geq x_0$, where $M>0$ is a constant. Equivalently, the Vinogradov notation $f\ll g$ means that $f(x)=O(g(x))$. We will use both these notations interchangeably as is convenient. 

We now discuss the sum
$S_z(x,y)$ (see \eqref{AMsum} below for definition) for complex $z$, before we state our results for $S_{-k}(x,y)$. But first, we need to state Selberg's classic estimate for $S_z(x)=S_z(x,x)$. 

In order to estimate the number of integers below a given magnitude with a prescribed number of prime factors, Selberg \cite{Sel} studied the sum
$$
S_z(x)=\sum_{1\leq n\leq x}z^{\nu(n)}
$$
for complex $z$. Since 
\begin{equation*}
\sum_{n=1}^\infty\dfrac{z^{\nu(n)}}{n^s}=\prod_{p}\left(1+\dfrac{z}{p^s-1}\right)=\zeta(s)^z\cdot f(s,z)    
\end{equation*}
where 
\begin{equation*}
f(s,z):=\prod_{p}\left(1+\dfrac{z}{p^s-1}\right)\left(1-\dfrac{1}{p^s}\right)^z    
\end{equation*}
is analytic for $\sigma=\operatorname{Re}(s)>1/2$, Selberg was able to show using the Perron integral method that 
\begin{theorem}[Selberg]\label{Sel}
For $R>0$, we have uniformly for $|z|\leq R$
\begin{equation*}
S_z(x)=\dfrac{x\log^{z-1}(x) f(1,z)}{\Gamma(z)}+O_R\left(x\log^{z-2}(x)\right).    
\end{equation*}
\end{theorem}
Using the above theorem, Selberg obtained uniform asymptotic estimates for the number $N_k(x)$ of integers $n\leq x$ with $k$ prime factors by expresssing $N_k(x)$ as a contour integral of $\frac{S_z(x)}{z^{k+1}}$, and choosing the contour suitably. 

Motivated by the work of the first author on $M(x,y)$ and $M^*(x,y)$, and Selberg's results on $S_z(x)$, Alladi and Molnar \cite{AM} studied the sum
\begin{align}\label{AMsum}
S_z(x,y):=\sum_{1\leq n\leq x}z^{\nu_y(n)}    
\end{align}
with the intent of estimating
$$
N_k(x,y)=\sum_{n\leq x, \nu_y(n)=k}1.
$$
Using the Perron integral method, the following result was established in \cite{AM}:

\begin{theorem}[Alladi-Molnar]\label{AM1}
Let $R>0$ be fixed and $|z|\leq R$ and $3\leq y\leq x$, then
\begin{eqnarray*}
S_z(x,y)=x\prod_{p<y}\left(1+\dfrac{z-1}{p}\right)+O\left(xe^{-\alpha}\log^D(x)\right)+O_R\left(\dfrac{x}{\log^{R+2}(x)}\right)
\end{eqnarray*}
for some absolute constant $D>0$. Thus, for $\epsilon>0$, if $\alpha\geq (R+D+1+\epsilon)\log\log x$ or $y\leq \exp\left(\dfrac{\log x}{(R+D+1+\epsilon)\log\log x}\right)$, then
\begin{eqnarray*}
S_z(x,y)=x\prod_{p<y}\left(1+\dfrac{z-1}{p}\right)+O\left(\dfrac{x}{\log^{R+1+\epsilon}(y)}\right)
\end{eqnarray*}
\end{theorem}
Thus if $z\neq 1-p$ for any prime $p$, we have
\begin{eqnarray*}
S_z(x,y)\asymp x\prod_{p<y}\left(1+\dfrac{z-1}{p}\right).
\end{eqnarray*}
\begin{remark}
Note that Theorem \ref{AM1} yields an asymptotic estimate for $S_z(x,y)$ only when
$$
z\neq 1-p \quad \text{for any} \quad p, \quad \text{and} \quad
y<x^{c/\log\log x}, \quad \text{for some} \quad c>0. 
$$
Thus Theorem \ref{AM1} yields an asymptotic estimate for $S_z(x,y)$ only for small $y$, where ``small'' is given by the above bound.
\end{remark}
Since $S_z(x,y)$ satisfies the Buchstab identity 
\begin{equation}\label{Buchrec}
S_z(x,y)=S_z(x,y^h)+(1-z)\sum_{y\leq p< y^h}S_z\left(\frac{x}{p},p\right),\hspace{1cm}h>1    
\end{equation}
Alladi and Molnar \cite{AM} were able to estimate $S_z(x,y)$ asymptotically
by induction on $\lfloor\alpha\rfloor$ by the Buchstab-de Bruijn method
for $y$ in a certain large range when $\operatorname{Re}(z)>0$ as enunciated in the
following theorem:
\begin{theorem}[Alladi-Molnar]\label{AM}
Let $\operatorname{Re}(z)>0$ and $2\leq \sqrt x\leq y\leq x$. Let $r=\operatorname{Re}(z)$. Then we have uniformly for $|z|\leq R$
\begin{eqnarray*}
S_z(x,y)=\dfrac{x\cdot m_z(\alpha)}{\log^{1-z} y}+O\left(\dfrac{(1-z)x}{\log^{2-r} y}+\dfrac{(1-z)x}{\log x}\right)
\end{eqnarray*}
where
\begin{eqnarray*}
m_z(\alpha)=\dfrac{g(1,z)}{\Gamma(z)}\left(\dfrac{1}{\alpha^{1-z}}+\dfrac{1-z}{\alpha^{1-z}}\int_1^\alpha\dfrac{du}{u^z(u-1)^{1-z}}\right).
\end{eqnarray*}
On the other hand, if $\operatorname{Re}(z)=r>0$ and $\alpha >2$ is fixed, then uniformly for $|z|\leq R$, we have
\begin{eqnarray*}
S_z(x,y)=\dfrac{x\cdot m_z(\alpha)}{\log^{1-z} y}+O\left(\dfrac{x\log\log x}{\log^{2-r} y}+\dfrac{x\log\log x}{\log y}\right)
\end{eqnarray*}
where 
\begin{eqnarray*}
m_z(\alpha)=\dfrac{2^{1-z}m_z(2)}{\alpha^{1-z}}+\dfrac{1-z}{\alpha^{1-z}}\int_1^\alpha\dfrac{m_z(u-1)\;du}{u^z}.
\end{eqnarray*}
\end{theorem}
We point out that a uniform version of Theorem \ref{AM} valid for
$$
\exp{((\log x)^{1-\delta})}\leq y\leq x
$$
with a certain $\delta >0$ is established in \cite{AM}. 
It turns out that the ranges of validity of Theorem \ref{AM1} and the uniform
version of Theorem \ref{AM} overlap and so the entire interval $2\leq y\leq x$ is covered. Using these theorems, and suitable bounds for $S_z(x,y)$ when $\operatorname{Re}(z)\leq 0$
that follow from the work of Tenenbaum \cite{T17}, the problem of the asymptotic estimation of $N_k(x,y)$ was resolved in \cite{AM} by the Selberg contour integral method. It is to be noted that the asymptotic behavior of $N_k(x,y)$ exhibits a certain
interesting variation of the classical theme as observed initially by Alladi, and this variation also prompted the study of $S_z(x,y)$ for complex $z$. Motivated by this variation of the classical theme, Tenenbaum \cite{T} used the
saddle point method to estimate $S_z(x,y)$ to get uniform asymptotic estimates for $N_k(x,y)$. 

It is to be noted that Selberg's theorem (Theorem \ref{Sel}) does not yield an asymptotic estimate
for $S_{-k}(x)$, when $k$ is a positive integer, because the main term vanishes
owing to the presence of $\Gamma(-k)$ in the denominator. Indeed
$S_{-k}(x)=O\left(x\cdot e^{-c\sqrt{\log\,x}}\right)$ with some constant $c>0$
follows from the Perron integral method and the zero-free region of the
Riemann zeta function. 
This makes the asymptotic study of $S_{-k}(x,y)$ extremely interesting. Since $S_{-k}(x,y)$ satisfies the Buchstab recurrence
\begin{align}\label{actBuchstab}
S_{-k}(x,y)=S_{-k}(x,y^h)+(k+1)\sum_{y\leq p<y^h}S_{-k}\left(\frac{x}{p},p\right)
\end{align}
the upper bound for $S_{-k}(x)=S_{-k}(x,x)$ can still be used to start the Buchstab iteration to estimate $S_{-k}(x,y)$
asymptotically. Now that we have provided the necessary background, we can describe our recent results on $S_{-k}(x,y)$, which we do in the next section.

\section{Statement of results on $S_{-k}(x,y)$}\label{Sect3} 

\noindent Here and in what follows, we shall use the following notation. We set
\begin{align*}
\beta=x/y,\quad \text{and} \quad \beta_\ell=x/y^\ell
\end{align*}
for an integer $\ell >1$. Also $R(t)$ will denote any decreasing function of $t$ that satisfies $R(t)=O\left(e^{-c\sqrt{\log\,t}}\right)$, with some
constant $c>0$ whose value does not concern us. Also, $c$ need not be the
same when used in diffrent contexts. Such an $R(t)$ can be used to bound
from above the relative error in the Prime Number Theorem, namely,
$$
\frac{|\pi(t)-\text{li}(t)|}{\text{li}(t)}<R(t).
$$
where $\pi(t)=\sum_{p\leq t}1$ and $\text{li}(t)=\int_2^t\frac{du}{\log u}$.

We begin by describing the asymptotic behavior of $S_{-k}(x,y)$ when
$\sqrt x <y\leq x$.  
In this region, the asymptotic
estimate for $S_{-k}(x,y)$ has several terms starting with an $x/\log y$ term, then an $x/\log^2 y$ term, and so on, up to an $x/\log^{k+1}y$ (resp. an $x/\log^{k+2}y$) term when $k\neq p-1$ (resp. $k=p-1$). Each of these terms (except the last term) is 
multiplied by a function of $\beta$, where these functions of $\beta$ rapidly
decay to zero, like the relative error $R(\beta)$ in the PNT as $\beta\rightarrow\infty$ with $x$. To be more precise, we have the following two theorems for $S_{-k}(x,y)$ for the cases $k\neq p-1$, for any prime $p$, and $k=p-1$ for some prime $p$:

\begin{theorem}\label{main1}
Let $k+1\neq p$ for any prime $p$. Then for $\sqrt{x}\leq y< x$ and $\beta=x/y$, we have
\begin{eqnarray*}
&&S_{-k}(x,y)=\dfrac{(k+1)xf_1(\beta)}{\log y}-\dfrac{(k+1)xf_2(\beta)}{\log^2y}+\dfrac{2(k+1)xf_3(\beta)}{\log^3 y}+\cdots\notag\\\mbox{}\\&&+\;(-1)^{k-1}\dfrac{(k+1)x(k-1)!f_k(\beta)}{\log^k y}+\dfrac{b_k(k+1)x}{\log^{k+1}x}+O\left(\dfrac{xR(\beta)}{\log^{k+1}y}\right)+O\left(\dfrac{x}{\log^{k+2}x}\right)
\end{eqnarray*}
where $b_k$ is the non-zero constant
\begin{equation*}
b_k=k!\prod_{p}\dfrac{p^k(p-k-1)}{(p-1)^{k+1}}  
\end{equation*}
and 
\begin{equation}
f_1(\beta):=\int_1^\beta\dfrac{S_{-k}(t)}{t^2}dt=O(R(\beta)),\hspace{1cm}f_j(\beta):=\int_1^\beta\dfrac{f_{j-1}(t)}{t}dt=O(R(\beta)),\hspace{0.8cm}2\leq j\leq k.    
\end{equation}
\end{theorem}
The important difference in the case $k=p-1$ for some $p=$ prime, is illustrated in the next result: 
\begin{theorem}\label{main2}
Let $k+1=p$ for some prime $p$. Then for $\sqrt{x}\leq y< x$ and $\beta=x/y$, we have
\begin{eqnarray*}
&&S_{-k}(x,y)=\dfrac{(k+1)xf_1(\beta)}{\log y}-\dfrac{(k+1)xf_2(\beta)}{\log^2y}+\dfrac{2(k+1)xf_3(\beta)}{\log^3 y}+\cdots\notag\\\mbox{}\\&&+\;(-1)^{k}\dfrac{(k+1)xk!f_{k+1}(\beta)}{\log^{k+1} y}+\dfrac{c_k(k+1)x}{\log^{k+2}x}+O\left(\dfrac{xR(\beta)}{\log^{k+2}y}\right)+O\left(\dfrac{x}{\log^{k+3}x}\right)
\end{eqnarray*}
where $c_k$ is the non-zero constant
\begin{eqnarray*}
c_k=(k+1)!\left(1+\dfrac{1}{k}\right)^{k+1}\log(k+1)\prod_{p\neq k+1}\dfrac{p^k(p-k-1)}{(p-1)^{k+1}}
\end{eqnarray*}
and
\begin{equation*}
f_1(\beta):=\int_1^\beta\dfrac{S_{-k}(t)}{t^2}dt=O(R(\beta)),\hspace{1cm}f_j(\beta):=\int_1^\beta\dfrac{f_{j-1}(t)}{t}dt=O(R(\beta)),\hspace{0.5cm}2\leq j\leq k+1.    
\end{equation*}
\end{theorem}
\begin{remark}
In Theorem 3.2, the functions $f_j(\beta)$ are given in the same manner, except that, we also have $f_j(\beta)=O(R(\beta))$ when $j=k+1$. Since $R(\beta)$ tends to zero rapidly as $\beta\to\infty$, what Theorem \ref{main1} and \ref{main2} say is that when $\beta$ in large, as it happens when $y=x^{1/\alpha}$ with $1+\delta\leq\alpha<2$, for some $\delta>0$ fixed, then all terms having
the $f_j(\beta)$ as factors become insignificant, and so the dominating term is of order of magnitude $x/\log^{k+1} y$ in Theorem \ref{main1} and $x/\log^{k+2}y$ in Theorem \ref{main2}. Henceforth, by {\bf{Case 1}} we shall mean that $k\neq p-1$ for any prime $p$, and by {\bf{Case 2}} we mean the $k=p-1$ for some
prime $p$. 
\end{remark}
With the asymptotic behavior of $S_{-k}(x,y)$ having been precisely determined
in the interval $\sqrt{x}<y\leq x$, the Buchstab identity \eqref{actBuchstab} can be used to estimate $S_{-k}(x,y)$ by induction on $\lfloor\alpha\rfloor$, that is in intervals $\ell<\alpha\le \ell+1$ in succession. But what happens is that in the interval $x^{1/3}<y<x^{1/2}$, there is no $x/\log y$ term and so the expansion starts with the $x/\log^2y$ term. Next in the interval $x^{1/4}<y<x^{1/3}$, there are no terms of size $x/\log y$ or $x/\log^2y$, thus the expansion starts with $x/\log^3y$, and so on until we have the stabilizing term $x/\log^{k+1}y$ when $k\neq p-1$, and $x/\log^{k+2}y$, when $k=p-1$. This  appealing phenomenon is described in the following theorems:
\begin{theorem}\label{main3}
Let $k+1\neq p$ for any prime $p$ and $\ell$ be an integer in $[2,k+1]$. Then there exists functions $E_{1,\ell}^{(\ell)}(t), E_{1,\ell+1}^{(\ell)}(t),\cdots, E_{1,k}^{(\ell)}(t)$ with $E_{1,j}^{(\ell)}(\beta_\ell)=O(R(\beta_\ell))$ $(\ell\leq j\leq k)$ such that in the range $\sqrt[\ell+1]{x}\leq y< \sqrt[\ell]{x}$, we have
\begin{eqnarray*}
S_{-k}(x,y)=\sum_{j=\ell}^k \dfrac{x\cdot E_{1,j}^{(\ell)}(\beta_\ell)}{\log^j x}+\dfrac{x\cdot m_{1,k}(\alpha)}{\log^{k+1}y}+O\left(\dfrac{xR(\beta_\ell)}{\log^{k+1}x}\right)+O\left(\dfrac{x}{\log^{k+2}y}\right)
\end{eqnarray*}
where for non-integral $\alpha\in[1,k+1]$, $m_{1,k}(\alpha)$ satisfies the
difference-differential equation
\begin{eqnarray}\label{diffeqm1k}
\alpha\cdot m_{1,k}'(\alpha)+(k+1)\cdot m_{1,k}(\alpha)=(k+1)\cdot m_{1,k}(\alpha-1),\hspace{0.5cm}\alpha>1
\end{eqnarray}
and more precisely, for $\ell=2,3,\cdots k$ and $\ell< \alpha< \ell+1$ the integral equation
\begin{eqnarray}\label{inteq1}
m_{1,k}(\alpha)=\dfrac{m_{1,k}(\ell^-)\ell^{k+1}}{\alpha^{k+1}}+d_{1,\ell}(k)\left(\dfrac{\ell}{\ell-1}\right)^k\dfrac{k+1}{\alpha^{k+1}}+\dfrac{k+1}{\alpha^{k+1}}\int_\ell^{\alpha}s^km_{1,k}(s-1)\;ds,
\end{eqnarray}
with a certain constant $d_{1,\ell}(k)$ specified below.
\end{theorem}

\begin{theorem}\label{main4}
Let $k+1=p$ for some prime $p$ and $\ell$ be an integer in $[2,k+2]$. Then there exists functions $E_{2,\ell}^{(\ell)}(t), E_{2,\ell+1}^{(\ell)}(t),\cdots, E_{2,k+1}^{(\ell)}(t)$ with $E_{2,j}^{(\ell)}(\beta_\ell)=O(R(\beta_\ell))$ $(\ell\leq j\leq k+1)$ such that in the range $\sqrt[\ell+1]{x}\leq y< \sqrt[\ell]{x}$, we have
\begin{eqnarray*}
S_{-k}(x,y)=\sum_{j=\ell}^{k+1} \dfrac{x\cdot E_{2,j}^{(\ell)}(\beta_\ell)}{\log^j x}+\dfrac{x\cdot m_{2,k}(\alpha)}{\log^{k+2}y}+O\left(\dfrac{xR(\beta_\ell)}{\log^{k+2}x}\right)+O\left(\dfrac{x}{\log^{k+3}y}\right)
\end{eqnarray*}
where for non-integral $\alpha\in[1,k+2]$, $m_{2,k}(\alpha)$ satisfies the
difference-differential equation
\begin{eqnarray}\label{diffeqm2k}
\alpha\cdot m_{2,k}'(\alpha)+(k+2)\cdot m_{2,k}(\alpha)=(k+1)\cdot m_{2,k}(\alpha-1),\hspace{0.5cm}\alpha>1
\end{eqnarray}
and more precisely, for $\ell=2,3,\cdots k+1$ and $\ell<\alpha<\ell+1$ the integral equation
\begin{eqnarray}\label{inteq2}
m_{2,k}(\alpha)=\dfrac{m_{2,k}(\ell^-)\ell^{k+2}}{\alpha^{k+2}}+d_{2,\ell}(k)\left(\dfrac{\ell}{\ell-1}\right)^{k+1}\dfrac{k+1}{\alpha^{k+2}}+\dfrac{k+1}{\alpha^{k+2}}\int_\ell^{\alpha}s^{k+1}m_{2,k}(s-1)\;ds,
\end{eqnarray}
with a certain constant $d_{2,\ell}(k)$ specified below.
\end{theorem}

It turns out that 
\begin{eqnarray*}
d_{1,\ell}(k):=\int_1^\infty \dfrac{E_{1,k}^{(\ell)}(t)}{t}dt\quad\text{and}\quad d_{2,\ell}(k):=\int_1^\infty \dfrac{E_{2,k+1}^{(\ell)}(t)}{t}dt. 
\end{eqnarray*}
Since the functions $E_{j,\ell}^{(\ell)}(\beta_\ell)$ ($j=1$ and $2$) in Theorems \ref{main3} and \ref{main4} tend rapidly to zero as $\beta_{\ell}\to\infty$, and in Theorems \ref{main1} and \ref{main2} the functions $f_j(\beta)\to 0$
rapidly as $\beta\to\infty$, it follows that if $\alpha$ is not an integer and fixed, then
\begin{eqnarray}
  S_{-k}(x,y)=\frac{x\cdot m_{1,k}(\alpha)}{\log^{k+1}y}+O\left(\frac{x}{\log^{k+2}y}\right), \quad
      \text{for} \quad \alpha\in (1,k+1), \quad \text{\bf{Case 1}},
\end{eqnarray}
and
\begin{eqnarray}
  S_{-k}(x,y)=\frac{x\cdot m_{2,k}(\alpha)}{\log^{k+2}y}+O\left(\frac{x}{\log^{k+3}y}\right), \quad
      \text{for} \quad \alpha\in (1,k+2), \quad \text{\bf{Case 2}}.
\end{eqnarray}

For {\bf{Case 1}}, when $\alpha > k+1$ (resp. for {\bf{Case 2}}, when
$\alpha >k+2$), all terms involving $x/\log^j y$ for $1\leq j\leq k$
(resp. $1\leq j\leq k+1$) disappear  when the Buchstab recurrence \eqref{actBuchstab} is
applied to estimate $S_{-k}(x,y)$ by induction on $\lfloor\alpha\rfloor$ (this is due to the properties of the function $G_k(s)$ discussed in the next section). Thus the the asymptotic behavior stabilizes. This is given by the following 
\begin{theorem}\label{main5}
(i) Let $\alpha>k+1$ be fixed. If $k\neq p-1$ for any prime $p$, then we have
\begin{equation*}
S_{-k}(x,y)=\dfrac{x\cdot m_{1,k}(\alpha)}{\log^{k+1}y}+O\left(\dfrac{x}{\log^{k+2}y}\right),
\end{equation*}
where $m_{1,k}(\alpha)$ satisfies the difference-differential equation \eqref{diffeqm1k}.

(ii) Let $\alpha>k+2$ be fixed. If $k=p-1$ for some prime $p$, then we have
\begin{equation*}
S_{-k}(x,y)=\dfrac{x\cdot m_{2,k}(\alpha)}{\log^{k+2}y}+O\left(\dfrac{x}{\log^{k+3}y}\right),
\end{equation*}
where $m_{2,k}(\alpha)$ satisfies the difference-differential equation \eqref{diffeqm2k}.
\end{theorem}

\begin{remark}
Although Theorem \ref{main5} is stated for fixed $\alpha$, it is possible to use
the Buchstab formula \eqref{actBuchstab} together with the integral equations satisfied by
$m_{1,k}(\alpha)$ and $m_{2,k}(\alpha)$, follow the method of Alladi in \cite{A1}, or of Alladi-Molnar in \cite{AM}, to show that for any $\varepsilon >0$, 
\begin{align}
\left|S_{-k}(x,y)-\dfrac{x\cdot m_{1,k}(\alpha)}{\log^{k+1}y}\right|\ll_{\varepsilon}
\dfrac{x\cdot {\alpha}^{k+2}}{\log^{k+2}y},\hspace{0.2cm} \text{for}\hspace{0.2cm} \alpha\geq k+1+\varepsilon \hspace{0.2cm} \text{(\bf{Case 1})},
\end{align}
and
\begin{align}
\left|S_{-k}(x,y)-\dfrac{x\cdot m_{2,k}(\alpha)}{\log^{k+2}y}\right|\ll_{\varepsilon}
\dfrac{x\cdot {\alpha}^{k+3}}{\log^{k+3}y},\hspace{0.2cm} \text{for} \hspace{0.2cm} \alpha\geq k+2+\varepsilon \hspace{0.2cm}  \text{(\bf{Case 2})}.
\end{align}

Thus we get an asymptotic estimate for $S_{-k}(x,y)$ when
$y>\exp{(\log^{1-\delta}(x))}$, for some $\delta >0$. However, to fully understand the above estimates, we
need to know the asymptotic behavior of the functions $m_{1,k}(\alpha)$ and
$m_{2,k}(\alpha)$ as $\alpha\to\infty$. For this purpose, and also to construct a continuous approximation to $S_{-k}(x,y)$, much like $\text{li}(x)$ is a continuous approximation to $\pi(x)$, we need to study certain Dirichlet series, which is what we do next. These Dirichlet series also shed light on the
constants $b_k$ and $c_k$ in Theorems \ref{main1} and \ref{main2} respectively. 
\end{remark}

\section{Some special Dirichlet series}

\noindent We start by defining a function $G_k(s)$ given as a convergent
Dirichlet series in $\operatorname{Re}(s)>1$, and which enjoys a second representation which is analytic in $Re(s)\ge 1$:
\begin{equation}\label{Eulerprod}
G_k(s):=\sum_{n=1}^\infty\dfrac{(-k)^{\nu(n)}}{n^s}=\prod_{p}\left(1-\dfrac{k}{p^s-1}\right)=\dfrac{g_k(s)}{\zeta^k(s)}    
\end{equation}
where 
\begin{equation*}
g_k(s):=\prod_p\dfrac{p^{sk}(p^s-(k+1))}{(p^s-1)^{k+1}}    
\end{equation*}
converges for $\operatorname{Re}(s)>1/2$, and $1/\zeta^k(s)$ is analytic in
$\operatorname{Re}(s)\ge 1$. Since the Dirichlet series is absolutely convergent for $Re(s)>1$ and uniformly convergent
in compact subsets of the half-plane $Re(s)>1$, we may differentiate term-by-term in that half-plane to deduce that
$$
G^{(j)}_k(s)=\sum^{\infty}_{n=1}\frac{(-k)^{\nu(n})(-1)^jlog^j(n)}{n^s}, \quad \operatorname{Re}(s)>1.
$$
However, each of Dirichlet series $G^{(j)}_k(s)$ for $j\ge 1$, converges at $s=1$. Indeed, by Steiltjes integration, we have
$$
G_k(1)=G^{(0)}_k(1)=\sum^N_{n=1}\frac{(-k)^{\nu(n)}}{n}
$$
$$
=\int^N_{1^-}\frac{dS_{-k}(x)}{x}
=\frac{S_{-k}(N)}{N}+\int^N_1\frac{S_k(x)dx}{x^2},
$$
because $S_{-k}(1^-)=0$. The bound $S_{-k}(x)=O(xR(x))$ shows that
$S_{-k}(N)/N$ tends to 0 as $N\to\infty$, and that
$$
\int^{\infty}_1\frac{S_{-k}(x)dx}{x^2} < \infty.
$$
It follows by the same argument with a little extra detail, that each of the series
$$
G^{(j)}_k(1)=\sum^{\infty}_{n=1}\frac{(-k)^{\nu(n)}(-1)^j\log^j(n)}{n},
$$
will converge for $j=1,2,\cdots$. The actual values of $G^{(j)}_k(1)$ are given below using the second representation; these values can be established by using ideas underlying Axer's theorem, 
the Perron integral method, and the zero-free region for $\zeta(s)$ in the
critical strip. 

We have two cases:

\textbf{Case\;1}:\; $k+1\neq p$ for any prime $p$.
\\

In this case, $G_k(s)$ has a zero of order $k$ at $s=1$ since $g_k(1)\neq 0$. Thus, we have
\begin{align*}
G_k(1)=G_k'(1)=G_k''(1)=\cdots=G_k^{(k-1)}(1)=0\hspace{0.5cm}\text{and}\hspace{0.5cm}G_k^{(k)}(1)\neq 0.    
\end{align*}
This means 
\begin{align*}
\displaystyle\sum_{n=1}^\infty\dfrac{(-k)^{\nu(n)}\log^j n}{n}=0,\hspace{0.5cm}
  \text{for}\quad j=0,1,2,\cdots,k-1, \quad \text{and} \quad \neq 0, \quad
  \text{if} \quad j=k. 
\end{align*}
  The value of $G_k^{(k)}(1)$ equals the constant $b_k$ in Theorem \ref{main1}, and we have
\begin{equation}\label{bk}
G_k^{(k)}(1)=b_k=k!\prod_{p}\dfrac{p^k(p-k-1)}{(p-1)^{k+1}}.    
\end{equation}

\textbf{Case\;2}:\; $k+1=p$ for some prime $p$.
\\

In this case, $G_k(s)$ has a zero of order $k+1$ at $s=1$ since $g_k(s)$ has a simple zero at $s=1$.
Thus, we have
\begin{align*}
G_k(1)=G_k'(1)=G_k''(1)=\cdots=G_k^{(k)}(1)=0\hspace{0.5cm}\text{and}\hspace{0.5cm}G_k^{(k+1)}(1)\neq 0.
\end{align*}
This means 
\begin{align*}
\displaystyle\sum_{n=1}^\infty\dfrac{(-k)^{\nu(n)}\log^j n}{n}=0,\hspace{0.6cm}
  \text{for}\quad j=0,1,2,\cdots, \quad \text{and}\quad \neq 0 \quad \text{if}\quad
  j=k+1.
\end{align*}

First, we consider \textbf{Case\;2} which is the exceptional case. We evaluate $G_k^{(k+1)}(1)$ in this case and show that this equals the constant $c_k$ in Theorem \ref{main2}. We have
\begin{align}\label{expconstant}
&G_k^{(k+1)}(1):=c_k=(k+1)!\left(1+\frac{1}{k}\right)^{k+1}\log(k+1)\prod_{p\neq k+1}\dfrac{p^k(p-k-1)}{(p-1)^{k+1}}.
\end{align}
Also, it can be shown that for any large $\beta\ge 1$
\begin{align}\label{fsum}
\sum_{n\leq \beta}\dfrac{(-k)^{\nu(n)}\log^{k+1} n}{n}=c_k+O\left(R(\beta)\right).
\end{align}
Indeed, using Steiltjes integration, \eqref{Eulerprod} and \eqref{expconstant}, we have
\begin{eqnarray}\label{steint}
\sum_{n\leq \beta}\dfrac{(-k)^{\nu(n)}\log^{k+1} n}{n}&=&\sum_{n=1}^\infty\dfrac{(-k)^{\nu(n)}\log^{k+1} n}{n}-\sum_{n>\beta}\dfrac{(-k)^{\nu(n)}\log^{k+1} n}{n}\notag\\&=&c_k-\int_{\beta^{-}}^\infty\dfrac{\log^{k+1} x}{x}d(S_{-k}(x))\notag\\&=&c_k-\left.\dfrac{\log^{k+1}(x)S_{-k}(x)}{x}\right|_{\beta^-}^\infty-\int_\beta^\infty\left(\dfrac{(k+1)\log^k x}{x^2}-\dfrac{\log^{k+1} x}{x^2}\right)S_{-k}(x)\;dx.
\end{eqnarray}
At this point, we use the estimate $S_{-k}(x)=O\left(x\cdot e^{-c\sqrt{\log\;x}}\right)$  to
deduce \eqref{fsum}. We also need the integrals $f_j(\beta)$ defined in Section \ref{Sect3}, as well as the
following integrals defined using the $f_j(\beta)$:
\begin{eqnarray*}
F_{j,\ell}(\beta):=\int_1^\beta\dfrac{f_j(t)\log^\ell t}{t}dt,\quad 1\leq j\leq k+i\quad \text{and}\quad 1\leq\ell\leq k-j+i.
\end{eqnarray*}
It turns out that for \textbf{Case\;2}, we have 
\begin{align}
&f_j(\beta)=O(R(\beta))\quad (1\leq j\leq k+1),\quad f_{k+2}(\beta)=\dfrac{(-1)^{k+1}c_k}{(k+1)!}+O(R(\beta)),\label{fj}\\
&\text{and}\hspace{0.7cm}F_{j,\ell}(\beta)=\begin{cases}
\dfrac{(-1)^{j}\ell!c_k}{(k+1)!}+O(R(\beta)),&\ell=k-j+1\\\mbox{}\\
O(R(\beta)),&\text{otherwise}.
\end{cases}\label{Fjl}
\end{align}
To deal with \textbf{Case\;1}, we consider similar integrals with $k+1$ replaced by $k$. 

\section{Sketch of the proof of Theorem \ref{main1}}
We only outline the proof of Theorem \ref{main1}. First, we use an an integral analogue of the Buchstab recurrence \eqref{actBuchstab} which follows from
the strong form of the Prime Number Theorem:
\begin{equation}\label{contbuchstab}
S_{-k}(x,y)=S_{-k}(x,y^h)+(k+1)\int_y^{y^h}S_{-k}(x/t,t)\dfrac{dt}{\log t}+O\left(x\cdot\alpha\cdot\log^k x\cdot R(y)\right).   
\end{equation}
Using \eqref{contbuchstab}, $S_{-k}(x,y)$ can be estimated first in the range $\sqrt{x}\leq y<x$ or, equivalently, $1<\alpha\leq 2$. To do so, we set $y^h=x$ in \eqref{contbuchstab} to obtain
\begin{equation}\label{setbuchstab}
S_{-k}(x,y)=S_{-k}(x,x)+(k+1)\int_y^{x}S_{-k}(x/t,t)\dfrac{dt}{\log t}+O\left(x\cdot\alpha\cdot\log^k x\cdot R(y)\right).    
\end{equation}
Noting that $S_{-k}(x,x)=S_{-k}(x)=O(x\cdot R(x))$, and $1\leq x/t\leq \sqrt{x}$ since $\sqrt{x}\leq y\leq t\leq x$, we immediately see that $S_{-k}(x/t,t)=S_{-k}(x/t)$. Thus, \eqref{setbuchstab} yields
\begin{equation}\label{buchsim}
S_{-k}(x,y)=(k+1)\int_y^{x}S_{-k}(x/t)\dfrac{dt}{\log t}+O(x\cdot R(y)).    
\end{equation}
By making the change of variable $t=x/u$, and applying integration by parts repeatedly by choosing suitable anti-derivatives determined by the functions $f_j(\beta)$ ($1\leq j\leq k$) in Theorem \ref{main1}, we end up with the following:
\begin{eqnarray}\label{intpartsrepeat}
\int_y^{x}S_{-k}(x/t)\dfrac{dt}{\log t}&=&\dfrac{x\cdot f_1(\beta)}{\log y}-\dfrac{x\cdot f_2(\beta)}{\log^2 y}+\cdots+(-1)^{k-1}(k-1)!\dfrac{x\cdot f_{k}(\beta)}{\log^{k} y}\notag\\&&+\;x(-1)^{k}k!\int_1^\beta\dfrac{f_{k}(s)}{s\log^{k+1} (x/s)}ds.  
\end{eqnarray}
Next, it follows by induction, and the fact that $S_{-k}(\beta)=O(\beta\cdot R(\beta))$ that $f_j(\beta)=O(R(\beta))$ for $1\leq j\leq k$. If we now set $f_{k+1}(\beta):=\int_1^\beta\dfrac{f_k(s)}{s}ds$, then these estimates coupled with repeated integration by parts yield 
\begin{eqnarray}\label{intp1}
f_{k+1}(s)&=&\dfrac{(-1)^{k-1}}{(k-1)!}\int_1^s\dfrac{f_1(t)\log^{k-1}t}{t}dt+O(R(s))\notag\\
&=&\dfrac{(-1)^{k-1}}{(k-1)!}\int_1^s f_1(t)\log^{k-1} td(\log t)+O(R(s))\notag\\
&=&\dfrac{(-1)^{k-1}}{(k-1)!}\left\{\left.\dfrac{f_1(t)\log^{k} t}{k}\right|_1^s - \dfrac{1}{k}\int_1^s\log^{k}t d(f_1(t))\right\}+O(R(s)).
\end{eqnarray}
The first expression inside curly braces in the right-hand side of \eqref{intp1} is $O(R(s))$. To evaluate the integral inside curly braces in the right-hand side of \eqref{intp1}, we need only observe that $d(f_1(t))=\frac{S_{-k}(t)}{t^2}dt$ whence changing the order of sum and integral, and repeated integration by parts yield
\begin{eqnarray}\label{logint}
\int_1^s\log^{k}t d(f_1(t))=\sum_{n\leq \beta}(-k)^{\nu(n)}\int_n^s\dfrac{\log^{k}t}{t^2} dt=b_k+O(R(s))
\end{eqnarray}
where $b_k$ is as in \eqref{bk}. Thus, \eqref{intp1} and \eqref{logint} yield
\begin{equation*}
f_{k+1}(s)=\dfrac{(-1)^k}{k!}b_k+O(R(s)),    
\end{equation*}
and this allows us to estimate the integral in the right-hand side of \eqref{intpartsrepeat} by choosing the following anti-derivative: 
\begin{equation}\label{finalintp}
\dfrac{f_{k}(s)ds}{s}=d\left(f_{k+1}(s)-\dfrac{(-1)^k}{k!}b_k\right)\quad \text{where}\quad f_{k+1}(s)-\dfrac{(-1)^k}{k!}b_k\rightarrow 0,\;s\rightarrow\infty.  
\end{equation}
Using \eqref{finalintp} in the right-hand side of \eqref{intpartsrepeat}, and applying integration by parts for the last time yields
\begin{eqnarray*}
\int_y^{x}S_{-k}(x/t)\dfrac{dt}{\log t}&=&\dfrac{x\cdot f_1(\beta)}{\log y}-\dfrac{x\cdot f_2(\beta)}{\log^2 y}+\cdots+(-1)^{k-1}(k-1)!\dfrac{x\cdot f_{k}(\beta)}{\log^{k} y}\notag\\&&+\;\dfrac{b_k\cdot x}{\log^{k+1}x}+O\left(\dfrac{x\cdot R(\beta)}{\log^{k+1}x}\right)+O\left(\dfrac{x}{\log^{k+2}y}\right)    
\end{eqnarray*}
which combines with \eqref{buchsim} to yield the required result, and we are done.

Theorem \ref{main2} follows analogously, the only difference is the fact that the expansion will now have an additional main term involving $\log^{k+2} x$ in the denominator. To obtain Theorems \ref{main3} and \ref{main4}, we once again follow the same strategy. Noting that $1<\frac{\log(x/t)}{\log t}\leq 2$, \eqref{contbuchstab} allows us to get an estimate for $S_{-k}(x,y)$ for $\sqrt[3]{x}\leq y<\sqrt{x}$, or equivalently, $2<\alpha\leq 3$ using the estimates obtained in Theorems \ref{main1} and \ref{main2} respectively. Finally, using induction on $\lfloor\alpha\rfloor$, and the Buchstab identity \eqref{contbuchstab} yield asymptotic estimates for $S_{-k}(x,y)$ for $\sqrt[\ell+1]{x}\leq y<\sqrt[\ell]{x}$.

\section{Continuous approximation to $S_{-k}(x,y)$}
Next, we construct continuous approximations to $S_{-k}(x,y)$ in both \textbf{Case 1} and \textbf{Case 2} which yield an asymptotic estimate for $S_{-k}(x,y)$ for a longer range of $y$. 
For this, define
\begin{equation*}
w_{k}(\alpha):=\begin{cases}
0,&\alpha\leq 1\\\mbox{}\\
\dfrac{k+1}{\alpha},&1<\alpha<2\\\mbox{}\\
\dfrac{k+1}{\alpha}+\dfrac{k+1}{\alpha}\displaystyle\int_{2}^\alpha w_{k}(s-1)\;ds,&\alpha>2.
\end{cases}    
\end{equation*}
Then we get
\begin{equation*}
m_{1,k}(\alpha)=\dfrac{(-1)^{k}b_k}{k!}w_{k}^{(k)}(\alpha),\hspace{1cm}m_{2,k}(\alpha)=\dfrac{(-1)^{k+1}c_k}{(k+1)!}w_{k}^{(k+1)}(\alpha).    
\end{equation*}
Next, define the following:
\begin{equation}\label{Acont}
A_{i,k}(x,y):=\dfrac{x}{\log y}\int_1^\infty w_{k}\left(\dfrac{\log x-\log t}{\log y}\right)\;df_1(t)   
\end{equation}
where $A_{i,k}(x,y)$ for $i=1$ (resp. $i=2$) corresponds to {\bf{Case 1}}
(resp. {\bf{Case 2}}). The form of the integral in (5.1) is the same for
both cases. The distinction in the notation for $A_{i,k}(x,y)$ is
due to the difference in the expansion in (5.2) below. It can be shown that 
\begin{equation*}
A_{i,k}(x,y)=A_{i,k}(x,y^h)+(k+1)\int_y^{y^h}A_{i,k}\left(\frac{x}{t},t\right)\dfrac{dt}{\log t}+O(xR(y))    
\end{equation*}
which is a continuous analogue of the Buchstab recurrence \eqref{actBuchstab}. 
\begin{theorem}
For $i=1$ \emph{(}$k\neq p-1$\emph{)} or $i=2$ \emph{(}$k=p-1$\emph{)}, we have
\begin{equation*}
S_{-k}(x,y)-A_{i,k}(x,y)\ll x\alpha^{k+i}\log^{k+i-2} xR(y).    
\end{equation*}
\end{theorem}
One can integrate $A_{i,k}(x,y)$ by parts successively. At every stage in the integration, the proper constant in the anti-derivative has to be chosen to ensure convergence. This leads to the series expansion
\begin{eqnarray}
A_{i,k}(x,y)&=&\dfrac{x(k+1)f_1(\beta)}{\log y}-\dfrac{x(k+1)f_2(\beta)}{\log^2 y}+\cdots+(-1)^{k+i}\dfrac{x(k+1)(k+i-2)!f_{k+i-1}(\beta)}{\log^{k+i-1} y}\notag\\\mbox{}\notag\\&&+\;\dfrac{x\cdot m_{i,k}(\alpha)}{\log^{k+i}y}+O\left(\dfrac{x}{\log^{k+i+1}y}\right),\hspace{1cm}\sqrt{x}<y<x,\hspace{1cm}i=1,\;2.
\end{eqnarray}
Note that the above expansion for $A_{i,k}(x,y)$ is identical to the series representation
for $S_{-k}(x,y)$ in Theorems \ref{main1} and \ref{main2}. 
A similar expansion for $A_{i,k}(x,y)$ starting from $\frac{x}{\log^\ell y}$ in the interval $x^{1/(\ell+1)}<y<x^{1/\ell}$ can be obtained and it will be identical to the series in Theorems \ref{main3} and \ref{main4}.
The continuous approximations $A_{1,k}(x,y)$ and $A_{2,k}(x,y)$ yield
\begin{theorem}\label{contapp}Let $\varepsilon>0$ be arbitrary, but small.
Then for 
for $\exp((\log x)^{1-\delta})<y<x$ for some $\delta>0$, we have 
\begin{eqnarray*}
S_{-k}(x,y)=\begin{cases}
\dfrac{x\cdot m_{1,k}(\alpha)}{\log^{k+1}y}+O_{\varepsilon}\left(\dfrac{x}{\log^{k+2}y}\right),&k\neq p-1\;\text{for any prime}\;p\;(\textbf{\emph{Case 1}}),\\\mbox{}\\
\dfrac{x\cdot m_{2,k}(\alpha)}{\log^{k+2}y}+O_{\varepsilon}\left(\dfrac{x}{\log^{k+3}y}\right),&k=p-1\;\text{for some prime}\;p\;(\textbf{\emph{Case 2}})
\end{cases}
\end{eqnarray*}
where $|\alpha-j|\geq \varepsilon$, for $j=1, 2, \cdots, k+i$ \emph{(}$i=1$ \emph{(}\textbf{\emph{Case 1}}\emph{)} and $i=2$ \emph{(}\textbf{\emph{Case 2}}\emph{)}\emph{)} as $m_{i,k}(\alpha)$ has jump discontinuities at these points. 
\end{theorem}

\section{Asymptotic behavior of $m_{1,k}(\alpha)$ and $m_{2,k}(\alpha)$ as
  $\alpha\to\infty$}\label{Asympbeh}

\noindent For a function $f(\alpha)$ which satisfies
 \begin{equation}\label{diffeqf}
  f'(\alpha)=\frac{\lambda}{\alpha}\int^{\alpha}_{\alpha -1}f'(t)\;dt,\quad \text{for}
  \quad \alpha >\alpha_0,
 \end{equation}
where $\lambda$ is a constant, it follows by iteration on $\lfloor\alpha\rfloor$ that
 \begin{equation}\label{fdecay}
f'(\alpha)=O(e^{-\alpha\log\alpha+O(\alpha)}).     
 \end{equation}
See Alladi-Molnar \cite{AM} for details of the method to get the upper bound for $f'(\alpha)$ in (6.2). Since $f'(\alpha)$ tends to zero so rapidly, by writing $f(\alpha)$ as an integral of $f'$, it follows that
\begin{equation}
f(\alpha)=c+O\left(e^{-\alpha\log\alpha+O(\alpha)}\right), \quad \alpha\to\infty,
\end{equation}
where $c$ is a constant. 

With regard to the functions $m_{1,k}(\alpha)$, \eqref{diffeqm1k}
implies that if $f(\alpha)=m_{1,k}(\alpha)$, then $f'(\alpha)=m'_{1,k}(\alpha)$ satisfies \eqref{diffeqf} with $\lambda=-(k+1)$. Thus by \eqref{fdecay} we see that $m_{1,k}(\alpha)$ tends to a limit $c$ as $\alpha\to\infty$. With regard to $m_{2,k}(\alpha)$,
there is a crucial difference in the asymptotic behavior. To understand this,
consider an anti-derivative $f(\alpha)$ of $m_{2,k}(\alpha)$. In this case, the difference-differential equation \eqref{diffeqm2k} for $m_{2,k}(\alpha)$ shows that with this anti-derivative $f(\alpha)$, the function $f'(\alpha)=m_{2,k}(\alpha)$ satisfies \eqref{diffeqf} with $\lambda=-(k+1)$, that is
$$
m_{2,k}(\alpha)=-\frac{k+1}{\alpha}\int^{\alpha}_{\alpha -1}m_{2,k}(t)dt.
$$
Thus by \eqref{fdecay} we see that $m_{2,k}(\alpha) \to 0$ as $\alpha\to\infty$.  

In order to compute the limit of $m_{1,k}(\alpha)$ as $\alpha\to\infty$, let us consider the behavior of $S_{-k}(x,y)$ for ``small $y$''. Indeed by the Alladi-Molnar theorem (Theorem \ref{AM1}) on $S_z(x,y)$ for ``small $y$'' we have the following:
\begin{equation*}
S_{-k}(x,y)
\begin{cases}
\sim \displaystyle x\prod_{p<y}\left(1-\frac{k+1}{p}\right)\asymp \frac{x}{\log^{k+1}y},&\textbf{Case 1}\\\mbox{}\\
=O_R\left(\dfrac{x}{\log^R y}\right),&\textbf{Case 2}.
\end{cases}    
\quad\quad\quad\text{for}\quad k+2\leq y\leq x^{c/\log\log x}
\end{equation*}
where $R$ can be chosen arbitrarily large. 
Notice that the large $y$ interval in Theorem \ref{contapp} overlaps with the small $y$ interval $k+2\leq y\leq x^{c/\log\log x}$. So, we have all of $k+2\leq y\leq x$ covered. So by comparing the two estimates in an interval common to regions for ``small $y$'' and ``large $y$'', we see that $m_{1,k}(\alpha)$ approaches a non-zero limit whereas $m_{2,k}(\alpha)$ decays to zero as $\alpha\rightarrow\infty$ with $x$. More precisely, we have:\\

\noindent\textbf{Case\;1}:
$$
\lim_{\alpha\rightarrow\infty}m_{1,k}(\alpha)=\ell(-k):=e^{-(k+1)\gamma}\prod_{p}\left(1-\dfrac{k+1}{p}\right)\left(1-\dfrac{1}{p}\right)^{-k-1}\neq 0,
$$
and so 
$$
m_{1,k}(\alpha)=\ell(-k)+O(e^{-\alpha\log\alpha+O(\alpha)})   
$$
where $\gamma$ is the Euler's constant.\\

\noindent \textbf{Case\;2}:
$$
\lim_{\alpha\rightarrow\infty}m_{2,k}(\alpha)=0
$$
and so
$$
m_{2,k}(\alpha)\ll e^{-\alpha\log\alpha+O(\alpha)}.
$$
In this case, $\ell(-k)=0$ as can be seen from the product.
\begin{remark}
  If we did not have estimates in a common region to make a comparison, then we would have to use other techniques to determine the asymptotic behavior of $m_{1,k}(\alpha)$
  when $\alpha\to\infty$. This would have involved representing $m_{1,k}(\alpha)$
  as a Laplace
integral, and evaluate this by the saddle point method, and/or consider the
adjoint of the function $m_{1,k}(\alpha)$.
\end{remark}

In closing, we point out that subsequent to our work, de le Breteche and Tenenbaum \cite{deBre-Ten} have recently studied more generally, the sums of oscillating functions over smooth numbers, namely, integers up to $x$ all of whose prime factors are
$\le y$, and applied their methods to get sharper forms of our results on
$S_{-k}(x,y)$.

\medskip

\section*{Acknowledgement}
\noindent The first author would like to thank Prof. M. Vidyasagar for the kind invitation to speak at the Subbarao Centenary Conference which is when these results were presented. The first author also thanks Kaneenika Sinha and the other
conference organizers for arrangements pertaining to his conference lecture. The second author was an institute postdoctoral fellow at IIT Gandhinagar during 2020-21 under the project IP/IITGN/MATH/AD/2122/15 which is when this paper, based on work done in 2017-18, was written. He sincerely thanks the institute for the support. Finally, we thank the referee for a very careful reading of the manuscript and for helpful suggestions. 

\end{document}